\documentclass[11pt]{article}
\usepackage{latexsym,amsfonts,amssymb,amsmath,amsthm, verbatim}
\usepackage[all]{xy}
\theoremstyle{plain}
        \newtheorem{theo}{Theorem}
        \newtheorem{defi}{Definition}
        
        \newtheorem*{coro}{Corollary}
        
\theoremstyle{remark}
        \newtheorem*{exam}{Example}
        \newtheorem*{note}{Note}

\def\quasidets{For more on the quasideterminant, consult the survey article \cite{GGRW:1}.}

\def\antipodes{Under these conditions, one may uniquely define an antipode $s$ as follows. For $p$ in the first graded piece, put $s(p)=-p$ as you are obliged to do for primitive elements. For higher graded pieces, the definition is forced upon you by the coalgebra grading and induction.}

\title{\large $NSym\hookrightarrow \mathcal Q_{\infty}$ is not a Hopf map}
\author{Aaron Lauve}

\begin{document}
\maketitle

\begin{abstract}
In this note, we show that there is no Hopf algebra structure on $\mathcal Q_{\infty}$, the algebra of pseudo-roots of noncommutative polynomials, which extends the one existing on $NSym$ (one of its famous subalgebras).\end{abstract}
 
%
\section{Introduction}
The algebra $\mathcal Q_n$ was introduced in 2001 by Gelfand, Retakh, and Wilson as a model for factoring noncommutative polynomials \cite{GelRetWil:2}. It is a graded, quadratic algebra with the remarkable property of remaining Koszul despite having several large free subalgebras (cf. \cite{SerWil:1} and \cite{Pio:1} for two independent proofs). Because of its natural origins, and because of the combinatorial nature of its generators and relations, it is interesting to ask whether or not it is a Hopf algebra. 

Here we answer a related question inspired by the famous square of combinatorial Hopf algebras
$$
\xymatrix{
NSym \ar@{->>}[d]_{} \ar@{^{(}->}[r]^{} & \mathfrak S Sym \ar@{->>}[d]^{} \\
Sym \ar@{^{(}->}[r]_{} & QSym
}
$$
studied in great detail by a great many authors (cf. \cite{AguSot:1, BerMykSot:1, Ges:1, GesReu:1, Haz:2, LodRon:1, MalReu:1} and the references therein). $NSym_n$, the algebra of noncommutative symmetric functions in $n$ variables, may be naturally identified as a subalgebra of $\mathcal Q_n$. As a first step towards answering the larger question, we show that this identification cannot be extended to a map of Hopf algebras.

The next section is critical to the argument. It summarizes results in the theory of noncommutative polynomials and points to why the map $\Phi : NSym\hookrightarrow Q_{\infty}$ chosen later is the natural one. Sections \ref{sec:Qn} and \ref{sec:NSym} introduce the algebras $\mathcal Q_{\infty}$ and $NSym$ respectively. The final section contains the calculations showing that $\Phi$ is not a Hopf algebra map.

\subsection*{Notation}
When convenient, we use the combinatorists' notation. In particular $[n] = \{1,2,\ldots,n\}$. Also, we will have occasion to denote the collection of subsets of $\{1,2\ldots,n\}$ of cardinality $d$ as $\binom{[n]}{d}$, and to write $A\in \binom{[n]}{d}$ to mean $A$ is a particular subset of $[n]$ of cardinality $d$. Finally, we write $\gamma\models n$ when $\gamma$ is a \emph{composition} of $n$ (any sequence of positive integers which sum to $n$). 

Throughout this note, $D$ is a noncommutative field with center $F\supseteq\mathbb Q$. 

\section{Factoring Noncommutative Polynomials}\label{sec:polys}
A good reference for the first part of this section is Lam's book \cite{Lam:1}. Let $D$ be a division ring, and let $f\in D[t]$ be a monic polynomial of degree $n$ in one variable over $D$. Because $D$ is noncommutative, there is not a unique way to ``evaluate'' $f$ at an element $x\in D$. 
\begin{exam} Over the quaternions, with $f(t)=i+jt=i+tj$ and $t\mapsto k$ we have
$$
i+j(k) = 2i \quad\neq\quad  0 = i+(k)j\,.
$$
\end{exam}

We agree to evaluate $f$ by always first writing it as a \emph{left-}polynomial
$$
f(t) = a_0 + a_1 t + \cdots + a_{n-1}t^{n-1} + t^n\,,
$$
and then plugging in $x$
$$
f(x) = a_0 + a_1 x + \cdots + x^n.
$$

\begin{theo} If $f\in D[t]$ is a monic polynomial of degree $n$, then $x$ is a \emph{root} of $f$ (that is, $f(x)=0$) if and only if there exists a polynomial $g$ of degree $n-1$ such that $f(t) = g(t)(t-x)$.
\end{theo}

\begin{note} Remember that you must expand the expression $g(t)(t-x)$ before evaluating. So if $g(t) = b_0 + b_1t + \cdots + b_{n-2}t^{n-2} + t^{n-1}$ then the theorem is really asserting the equality of $f$ and 
$(b_0t + b_1t^2 + \cdots + t^{n}) - (b_0x + b_1xt + \cdots + b_{n-2}xt^{n-2} + xt^{n-1}).$
\end{note}

Expanding on this note, if a polynomial has a factorization $f(t) = g(t)h(t)$ it is generally not the case that $f(x) = g(x)h(x)$. In particular, roots of $g$ are not necessarily roots of $f$. 

\begin{exam} Over the quaternions, the polynomial $f(t) = t^2-(i+j)t - k$ has a factorization $f(t) = (t-j)(t-i)$ but exactly one root, $x=i$. \end{exam}

\begin{theo} Let $f(t)=g(t) h(t)$ be a factorization of $f$, and suppose $x\in D$ satisfies $h(x) = a$. If $a=0$, then $f(x) = g(x)h(x)$, otherwise, $f(x) = g(axa^{-1})h(x)$. \end{theo}

\begin{defi} We call the elements $y_r$ of $D$ showing up in a factorization $(t-y_n)(t-y_{n-1})\cdots(t-y_2)(t-y_1)$ of $f$ the \emph{pseudo-roots} of $f$, $y_1$ additionally being an actual root. \end{defi}

In \cite{GelRet:4} Gelfand and Retakh found a closed-form expression for the pseudo-roots involving the Vandermonde quasideterminant\footnote{\quasidets}.

\begin{defi} Given elements $x_1, x_2, \ldots, x_r$ in a division ring $D$, the \emph{Vandermonde} quasideterminant $V(x_{1},x_{2},\ldots,x_{r})$ is the $(1,r)$-th quasideterminant of the Vandermonde matrix:
$$
\left|\begin{array}{cccc}
x_{1}^{r-1} & x_{2}^{r-1} & \cdots & \framebox[1.1\width]{$x_{r}^{r-1}$}\\
\vdots &\vdots&& \vdots\\
x_{1}^{1} & x_{2}^{1} & \cdots & x_{r}^{1}\\
1 & 1 & \cdots & 1
\end{array}\right|.
$$
\end{defi}

\begin{defi}[Gelfand-Retakh, \cite{GelRet:4}] If $f [t]$ is a monic polynomial of degree $n$ with $n$ roots $x_1,\ldots,x_n$ then we say that the roots are \emph{independent} if $V(x_{i_1},x_{i_2},\ldots,x_{i_r})$ is defined for all $1\leq r\leq n$ and all orderings $(i_1,\ldots i_n)$ of the roots.
\end{defi}

\begin{theo}[Gelfand-Retakh, \cite{GelRet:4}] If $f$ has independent roots $x_1, \ldots, x_n$, then the pseudo-roots are given by the formulas:
\begin{eqnarray*}
y_1 &=& x_1 \\
y_2 &=& V(x_1,x_2)x_2V(x_1,x_2)^{-1}\\
&\vdots&\\
y_n &=& V(x_1,x_2,\ldots,x_n)x_nV(x_1,x_2\ldots,x_n)^{-1}.\\
\end{eqnarray*}
\end{theo}

In their paper, they go on to prove a noncommutative Vieta theorem.

\begin{theo} The following rational functions in the variables $x_1,\ldots, x_n$ are symmetric with respect to the $\mathfrak S_n$ action on the $x$ indices:
\begin{eqnarray}
e_1 &=& y_n + y_{n-1} + \cdots + y_1\label{eq:e1}\\
&\vdots&\\
e_r &=& \sum_{i_r>\cdots >i_1} y_{i_r}\cdots y_{i_2}y_{i_1}\\
&\vdots&\\
e_n &=& y_n y_{n-1}\cdots y_2 y_1\label{eq:en} .
\end{eqnarray}
\end{theo}
They conjecture that these are truly the elementary noncommutative symmetric functions. In other words, if $g$ is a rational function in the $x_i$ which is symmetric, then it is a polynomial in the $e_r$. 
This important conjecture is proven by Wilson in \cite{Wil:1}. Equally important to this note, he proves that the $y_r$, and hence  the $e_r$, are algebraically independent.

\section{The algebra $\mathcal Q_{\infty}$}\label{sec:Qn}
When the roots $x_1, \ldots x_n$ above are independent, we may change the order of the roots and get another (full) set of pseudo-roots. In particular, we should replace the symbol $y_r$ by something like $y_{(i_1,\ldots,i_n),r}$\ldots reserving $y_1, \ldots, y_n$ for the fixed ordering $(1,2,\ldots, n)$.

\begin{theo}[G-R-W, \cite{GelRetWil:2}] The symbols $y_{(i_1,\ldots,i_n),r}$ do not depend on the last $n-r$ roots, and do not depend on the ordering of the first $r-1$ roots. \end{theo}
So we settle on the notation $y_{\{i_1, \ldots, i_{r-1}\},i_r}$ for the collection of all pseudo-roots associated to a polynomial $f$.

Unlike the set $\{y_r \mid 1\leq r\leq n\}$, the set $\{y_{A,i} \mid A\in\binom{[n]}{r-1},\,i\in [n]\setminus A\}$ is not algebraically independent (not even linearly independent over $\mathbb Q$). In \cite{GelRetWil:2}, Gelfand, Retakh, and Wilson introduce the algebra $\mathcal Q_n$ as a model for the relationships between the pseudo-roots of polynomials of degree $n$. 

\begin{defi} Let $\mathcal Q_n$ be the algebra over $\mathbb Q$ with generators $\{x_{A,i} \mid A\in\binom{[n]}{r-1},\,1\leq r\leq n,\,i\in [n]\setminus A\}$ and relations
\begin{eqnarray}
\label{eq:sum} x_{A\cup i,j} + x_{A,i} & = & x_{A\cup j,i} + x_{A,j};\\
\label{eq:product} x_{A\cup i,j}\cdot x_{A,i} & = & x_{A\cup j,i}\cdot x_{A,j}.
\end{eqnarray}
for $i\neq j$ and $i,j\not\in A$. Let $\mathcal Q_{\infty}$ denote the direct limit of these algebras, the algebra with generators $\{x_{A,i} \mid A\subseteq\mathbb N, \,i\in \mathbb N\setminus A\}$ and relations given by (\ref{eq:sum}) and (\ref{eq:product}).
\end{defi}

Note that for any fixed ordering $(i_1,\ldots,i_n)$ of $[n]$, the subalgebra generated by $x_{\emptyset,i_1}, x_{\{i_1\},i_2}, \ldots, x_{\{i_1,i_2,\ldots,i_{n-1}\},n}$ is relation-free, as Wilson's theorem dictates. Note also that equation (\ref{eq:sum}), plus induction, allows us to throw away generators of the form $x_{A,i}$ when $\max{A}\not< i$. Eliminating these generators and relation (\ref{eq:sum}), we see that $Q_n$ is a quadratic algebra.

We introduce some notation to simplify the computations to come. Let $y_r$ be shorthand for the generator $x_{\{1,2,\ldots, r-1\},r}$, and let $e_1,\ldots,e_n$ denote the elements of $\mathcal Q_n$ given by equations (\ref{eq:e1}) through (\ref{eq:en}). Finally, when $A=\{a_1<a_2<\cdots<a_r\}$, let $X(A) = x_{\emptyset,a_1} + x_{\{a_1\},a_2} + \cdots + x_{\{a_1,a_2,\ldots,a_{r-1}\},a_r}$. In particular, $X([n]) = e_1$; put $X(\emptyset)=0$.

For the argument in Section \ref{sec:embedding}, we will need to know a bit more about $\mathcal Q_n$.

\subsection*{A basis for $\mathcal Q_n$:}
\begin{enumerate} 
\item The symbols $\{X(A) \mid A\subseteq[n]\}$ generate $\mathcal Q_n$. 

\item Suppose $A\in\binom{[n]}{r}$ and $0\leq j\leq r$. We define $A^{(j)} := \{a_1, a_2, \ldots, a_{r-j}\}$; i.e. $A$ with its last $j$ entries deleted.

\item Given $A\in\binom{[n]}{r}$ and $1\leq j\leq r$, write $(A:j)$ for the sequence  $(A^{(0)}, A^{(1)},\ldots,A^{(j-1)})$.

\item Suppose $B_i=(A_i:j_i),\,1\leq i\leq s$ be a collection of sequences of this type. Let $\mathcal=(B_1,\ldots,B_s)$ denote the concatenation of these sequences, $\mathcal B=(A_1,\ldots,A_1^{(j_1-1)},A_2,\ldots)$. Let us call such a concatenation a \emph{string}, and let $\mathrm{wt}\,\mathcal B = j_1+\cdots+j_s$ and $\ell(\mathcal B)=s$.

\item Writing $X(\mathcal B)$ for the product $X(A_1)\cdots X(A_1^{(j_1-1)})\cdots X(A_s^{(j_s-1)})$ we have:
\end{enumerate}

\begin{theo}[G-R-W, \cite{GelRetWil:2}] 
The set of all symbols $X(\mathcal B)$, $\ell(\mathcal B)=s$, where $j_i\leq |A_i|$ for all $i$, and for all $2\leq i\leq s$, we have either $|A_i|\neq |A_{i-1}|-j_{i-1}$ or $A_i\nsubseteq A_{i-1}$ is a basis for $Q_n$. 
\end{theo}

In other words, $X(\mathcal B)$ is not an element of our basis if and only if there exists $2\leq i \leq s$ such that $A_i\subseteq A_{i-1}$ and $|A_i|=|A_{i-1}|-j_{i-1}$.
More important for this note, we have the

\begin{coro} $Q_n$ is a graded, quadratic algebra with symbols $X(A)$ all having degree one, and the $i^{\mathrm{th}}$ graded piece spanned by those symbols $X(\mathcal B)$ allowed above satisfying $\mathrm{wt}\,\mathcal B = i$.
\end{coro}

Note that $y_r = X([r])-X([r-1])$ is homogeneous of degree one. We'll denote the $i^{\mathrm{th}}$ graded piece of $\mathcal Q_n$ by $\mathcal Q_{n,i}$.  Critical to the argument in Section \ref{sec:embedding} is the fact that 
$$
\mathcal Q_n \otimes \mathcal Q_n = \bigoplus_{(i,j)\in\mathbb N^2} \mathcal Q_{n,i} \otimes \mathcal Q_{n,j}\,.
$$

\section{The algebra $NSym$}\label{sec:NSym}
The algebra $NSym$ (over $\mathbb Q$) has an interesting history. It is isomorphic \cite{Haz:1} to one of the earliest Hopf algebras ever studied: the universal enveloping algebra of the free Lie algebra with countably many generators \cite{Reu:1}. It is isomorphic \cite{MalReu:1} to the sum of Solomon's descent algebras \cite{Sol:1} of type A. Neither of these explain its name. For that, let us recall the commutative algebra of symmetric functions $Sym$.

\begin{defi}\label{def:sym}$Sym$ is the collection of all functions $f$ on $\mathbb N$ variables which may be written as a polynomial in the elementary symmetric functions:
\begin{eqnarray*}
\Lambda_1 &=& x_1 + x_2 + \cdots \\
\Lambda_2 &=& x_1x_2 + x_1x_3 +\cdots +  x_2x_3 + \cdots \\
&\vdots&\\
\Lambda_r &=&\sum_{i_1<\cdots<i_r}x_{i_1}\cdots x_{i_r}\\
&\vdots&
\end{eqnarray*}
\end{defi}
This commutative $\mathbb Q$-algebra, freely generated by the $\Lambda_r$, has a bialgebra structure given by putting $\Delta(\Lambda_r) = \sum_{i+j=r} \Lambda_i\otimes\Lambda_j$ (where $\Lambda_0 = 1$). As it is a connected, graded bialgebra, it is automatically a Hopf algebra\footnote{\antipodes}. In this note, we show that $\Phi: NSym\hookrightarrow \mathcal Q_{\infty}$ is not a bialgebra map, so further discussion of antipodes will be omitted.

With Section \ref{sec:polys} and the discussion after definition \ref{def:sym} in mind, the following definition is now motivated.

\begin{defi} The algebra of noncommutative symmetric functions $NSym$ is the free noncommutative $\mathbb Q$-algebra with generators $z_i$ indexed by $\mathbb N$ ($z_0$ being identified with the unit in $NSym$). It is a connected, graded Hopf algebra with coalgebra structure $(\Delta,\varepsilon)$ given by $\Delta(z_r) = \sum_{i+j=r} z_i\otimes z_j$ and $\varepsilon(z_r) = \delta_{0r}$.
\end{defi}

In \cite{GKLLRT:1}, the authors develop the theory of $NSym$ in parallel to the theory of $Sym$. They use the quasideterminant to prove analogs of assorted theorems for $Sym$ which have determinantal proofs. Similar to the commutative version it has several important bases (complete and monomial symmetric functions, Schur functions, etc). The study of the structure constants for (co)multiplication with respect to these bases is a growing industry \cite{AguSot:1, LodRon:1, Rea:1}.

\section{The Embedding}\label{sec:embedding}
Let $NSym(n)$ denote the Hopf subalgebra generated by $\{z_0, \ldots, z_n\}$. The algebra $NSym$ is a direct limit of the algebras $NSym(n)$. In what follows, we concentrate not on the map $\Phi : NSym\hookrightarrow \mathcal Q_{\infty}$, but on the map $\Phi_n : NSym(n)\hookrightarrow \mathcal Q_{n}$ for a fixed $n$.

As defined, $NSym(n)$ is a free algebra on generators $z_1,\ldots,z_n$ , which we may call the ``noncommutative symmetric functions on $n$ variables.'' As we mentioned above, $e_1,\ldots,e_n$ \textbf{are} noncommutative symmetric functions on $n$ variables. We let $\Phi_n(z_r) = e_r\in\mathcal Q_n$, and show that this is not a bialgebra map. 

First, note that $\Phi$ is an injective algebra map since the $e_r$ are algebraically independent in $\mathcal Q_n$. Second, note that the image $\Phi(NSym)$ is contained in the subalgebra of $\mathcal Q_n$ generated by the pseudo-roots $y_1, \ldots, y_n$. Call this subalgebra $\mathcal P(n)$. 

\begin{theo} There is no bialgebra structure on $\mathcal Q_n$ which extends the bialgebra structure on $NSym(n)$.\end{theo}

\begin{coro} There is no Hopf algebra structure on $\mathcal Q_{\infty}$ which extends the Hopf algebra structure on $NSym$.\end{coro}

\begin{proof} Let $(\Delta, \varepsilon)$ be the coalgebra structure on $NSym_n$. We assume the existence of a coalgebra structure $(\tilde{\Delta}, \tilde{\varepsilon})$ on $\mathcal Q_n$ making $\Phi$ a coalgebra map. We needn't look past the generators of $\mathcal P(n)$ to reach a contradiction. Calculations will be presented for the case $n=3$ to make the exposition palatable. 

\subsection*{The counit map on the generators $y_r$:}
\begin{enumerate}
\item $\varepsilon(z_3) = 0$ and $\Phi(z_3)=e_3=y_3y_2y_1$ implies at least one of the $y_r$ must satisfy $\tilde{\varepsilon}(y_r)=0$ ($\mathcal P(3)$ is free algebra, in particular, a domain), say it's $y_3$.

\item $\varepsilon(z_2) = 0$ and $\Phi(z_2)=e_2=y_3y_2 + y_3y_1+y_2y_1$. After the assumption on $y_3$ above, we are left with $\tilde{\varepsilon}(y_2y_1)=0$. Suppose $\tilde{\varepsilon}(y_2)=0$.

\item $\tilde{\varepsilon}(z_1)=0$ and $\Phi(z_1)=y_3+y_2+y_1$ implies that $y_1$ is killed by $\tilde{\varepsilon}$ too.

\item Finally, relations (\ref{eq:sum}) and (\ref{eq:product}) imply that $\tilde{\varepsilon}(x_{A,i}) = 0$ for all generators $x_{A,i}$ of $\mathcal Q_n$
\end{enumerate}

\subsection*{The comultiplication map on the generators $y_r$:}
\begin{enumerate}
\item We begin in complete generality, putting 
$$\tilde{\Delta}(y_r)=\sum_{s\geq 0}\sum_{i+j=s}\sum_{{|\mathcal B|=i,}\atop{|\mathcal B'|=j}} C(r)_{\mathcal B,\mathcal B'} X({\mathcal B})\otimes X({\mathcal B'}).$$

\item Then, by the grading on $\mathcal Q_3\otimes\mathcal Q_3$,
$(\tilde{\varepsilon}\otimes 1)\tilde{\Delta}(y_r) = y_r=(1\otimes \tilde{\varepsilon})\tilde{\Delta}(y_r)$ implies 

\begin{enumerate}
\item $C(r)_{\emptyset,\mathcal B'}=C_{\mathcal B,\emptyset} = 0$ when 
$|\mathcal B|,|\mathcal B'|\neq 1$, 
\item $\sum_{|\mathcal B|=1} C_{\mathcal B,\emptyset}X(\mathcal B) =  y_r$,
\item $\sum_{|\mathcal B'|=1} C(r)_{\emptyset,\mathcal B'}X(\mathcal B') = y_r$,
\end{enumerate}

\item \emph{Conclude:} $\tilde{\Delta}(y_r) = 1\otimes y_r + y_r\otimes 1 + f_r,$ where $f_r$ is a linear combination of symbols $X({\mathcal B})\otimes X({\mathcal B'})$ belonging to $\bigoplus_{(i,j)\geq (1,1)} \mathcal Q_{3,i} \otimes \mathcal Q_{3,j}$.
\end{enumerate}
\smallskip

\noindent More explicitly, I'm claiming $X(\{1,2\})\otimes y_1$ may appear in $\tilde{\Delta}(y_2)$, but terms like $1\otimes 1$, $1\otimes X(\{1\})$, and $y_3y_1\otimes 1$ won't.
\medskip

We will show that there is no definition of $\tilde{\Delta}$ that satisfies $\Phi\circ\Delta z_2 = \tilde{\Delta}\circ\Phi(z_2)= \tilde{\Delta}e_2$. First we'll need to compute $ \tilde{\Delta}e_1.$ We know that $\Delta z_1 = z_1\otimes 1 + 1\otimes z_1$, so we must have
\begin{eqnarray*}
\tilde{\Delta}(e_1) & = & e_1\otimes 1 + 1\otimes e_1\\
& = & (y_3+y_2+y_1)\otimes 1 + 1\otimes(y_3+y_2+y_1),\\
\tilde{\Delta}(y_3+y_2+y_1)  & = & (1\otimes y_3 + y_3\otimes 1 + f_3) + (1\otimes y_2 + y_2\otimes 1 + f_2) \\
&&\quad +\, (1\otimes y_1 + y_1\otimes 1 + f_1).
\end{eqnarray*}
This implies 
\begin{equation}\label{eq:delta1}
f_1 + f_2 + f_3=0,
\end{equation}
an equation taking place in $\bigoplus_{(i,j)\geq (1,1)}\mathcal Q_{3,i}\otimes\mathcal Q_{3,j}$.

Now we consider $\tilde{\Delta}(e_2)$. As the image of $\Phi(\Delta z_2)$, it must satisfy
\begin{eqnarray*}
 1\otimes e_2 + e_1\otimes e_1 + e_2\otimes 1 & = & \tilde{\Delta}(e_2)\\
 & = & \tilde{\Delta}(y_3y_2+y_3y_1+y_2y_1)\\
 & = & (1\otimes y_3 + y_3\otimes 1 + f_3)(1\otimes y_2 + y_2\otimes 1 + f_2)\\
 && +(1\otimes y_3 + y_3\otimes 1 + f_3)(1\otimes y_1 + y_1\otimes 1 + f_1)\\
 && +(1\otimes y_2 + y_2\otimes 1 + f_2)(1\otimes y_1 + y_1\otimes 1 + f_1)\\
 & = & 1\otimes e_2 + e_1\otimes e_1 + e_2\otimes 1 + (\large{\star}),
 \end{eqnarray*}
where $(\large{\star})$, appearing below, must be zero.
\begin{equation}\label{eq:delta2}
f_3\tilde{\Delta}y_2+\tilde{\Delta}y_3f_2+f_3\tilde{\Delta}y_1+\tilde{\Delta}y_3f_1+f_2\tilde{\Delta}y_1+\tilde{\Delta}y_2f_1 -y_3\otimes y_3-y_2\otimes y_2-y_1\otimes y_1= 0.
\end{equation}

Using equation (\ref{eq:delta1}) to simplify equation (\ref{eq:delta2}) we get
\begin{eqnarray*}
(-f_1)(y_1\otimes 1 + 1\otimes y_1 + f_1) &&\\
+\, (y_3\otimes 1 + 1\otimes y_3 + f_3)(-f_3) &&\\
+\, (f_3)(y_2\otimes 1 + 1\otimes y_2 + f_2) &&\\
+\, (y_2\otimes 1 + 1\otimes y_2 + f_2)(f_1) &=& y_3\otimes y_3 + y_2\otimes y_2 + y_1\otimes y_1.
\end{eqnarray*}
Now, each term on the left belongs to $\bigoplus_{(i,j)>(1,1)}Q_{3,i}\otimes Q_{3,j}$ while each term on the right belongs to $Q_{3,1}\otimes Q_{3,1}$. So both must be zero. But the expression on the right (writing $X_i$ for $X([i])$) is equal to 
\begin{equation*}
X_3\otimes \left(X_3 - X_2\right)
 + X_2\otimes \left(2X_2 - X_3- X_1\right)  + X_1\otimes \left(2X_1-X_2\right),
\end{equation*}
which is clearly nonzero. A contradiction.
\end{proof}

\bibliographystyle{amsplain}   
\def\cprime{$'$} \def\cprime{$'$}
\providecommand{\bysame}{\leavevmode\hbox to3em{\hrulefill}\thinspace}
\providecommand{\MR}{\relax\ifhmode\unskip\space\fi MR }
\providecommand{\MRhref}[2]{%
  \href{http://www.ams.org/mathscinet-getitem?mr=#1}{#2}
}
\providecommand{\href}[2]{#2}

\end{document}